\documentclass[12pt]{article}
\usepackage{fullpage,graphicx,psfrag,amsmath,amsfonts,verbatim}
\usepackage[small,bf]{caption}\usepackage{amssymb}
\usepackage{enumitem}
\usepackage{tikz}
\usepackage{listings}
\usepackage{authblk}
\usepackage{textcomp}
\usepackage{makecell}
\graphicspath{{figures/}}

\definecolor{codegreen}{rgb}{0,0.6,0}
\definecolor{codegray}{rgb}{0.5,0.5,0.5}
\definecolor{codepurple}{rgb}{0.58,0,0.82}
\definecolor{backcolour}{rgb}{0.95,0.95,0.98}

\lstdefinestyle{mystyle}{
backgroundcolor=\color{backcolour},   
keywordstyle=\color{magenta},
numberstyle=\tiny\color{codegray},
stringstyle=\color{codepurple},
basicstyle=\ttfamily\small,
breakatwhitespace=false,         
breaklines=true,                 
captionpos=t,                    
keepspaces=true,                 
numbers=left,                    
numbersep=5pt,                  
showspaces=false,                
showstringspaces=false,
showtabs=false,                  
tabsize=2,
}

\definecolor{seagreen}{rgb}{0.18, 0.55, 0.34}
\definecolor{mediumviolet-red}{rgb}{0.78, 0.08, 0.52}
\definecolor{khaki}{rgb}{0.94, 0.9, 0.55}

\usepackage{hyperref}
\usepackage{wasysym}
\hypersetup{
colorlinks   = true,    
urlcolor     = blue,    
linkcolor    = blue,    
citecolor    = red      
}

\newcommand{\eg}{{\it e.g.}}
\newcommand{\ie}{{\it i.e.}}

\newcommand{\BEAS}{\begin{eqnarray*}}
\newcommand{\EEAS}{\end{eqnarray*}}
\newcommand{\BEA}{\begin{eqnarray}}
\newcommand{\EEA}{\end{eqnarray}}
\newcommand{\BEQ}{\begin{equation}}
\newcommand{\EEQ}{\end{equation}}
\newcommand{\BIT}{\begin{itemize}}
\newcommand{\EIT}{\end{itemize}}

\newcounter{algorithmctr}
\renewcommand{\thealgorithmctr}{\arabic{algorithmctr}}
   {\mbox{}\\*[\parskip]\begin{minipage}{\linewidth}%
       \refstepcounter{algorithmctr}\begin{list}{}{%
       \setlength{\rightmargin}{0\linewidth}%
       \setlength{\leftmargin}{.05\linewidth}}%
       \rmfamily\small
       \item[]{\setlength{\parskip}{0ex}\hrulefill\par%
        \nopagebreak{\bfseries\textsf{Algorithm \thealgorithmctr~}}}}%
   {{\setlength{\parskip}{-1ex}\nopagebreak\par\hrulefill\\*[2ex]\par}%
   \end{list}\end{minipage}}

\title{Aging-Aware Battery Control via Convex Optimization}

\author[1]{Obidike Nnorom Jr.}
\author[1]{Giray Ogut}
\author[1]{Stephen Boyd}
\author[2] {Philip Levis}
\affil[1]{Department of Electrical Engineering, Stanford University}
\affil[2]{Department of Computer Science, Stanford University}

\begin{document}
\maketitle
\begin{abstract}
We consider the task of controlling a battery while balancing two competing objectives
that evolve over different time scales.
The primary objective,  
such as generating revenue by exploiting time varying energy prices or 
smoothing out the load of a computation center, 
operates on the scale of hours or days.
The long term objective is to maximize the lifetime of the battery,
which operates on a time scale of months and years.
These objectives conflict; roughly speaking, the primary objective
improves with cycling the battery more, which ages the battery faster.
Using an existing model for battery aging,
we formulate the problem of controlling the battery
under these competing objectives as a convex optimization problem. 
We demonstrate the tradeoff
between the primary objective and battery lifetime through 
numerical simulations.
\end{abstract}
\newpage
\tableofcontents
\newpage

\section{Introduction}
\subsection{Setting and tasks}\label{ss-setting}
In the last decade, battery storage, whether retail or grid scale, has become
increasingly common. The ability to store energy and discharge it at a later
time allows for temporal decoupling of energy production and consumption.
Providing this service can be lucrative, as the price of electricity can vary by 
orders of magnitude depending on the time of day and the season. 
In other basic applications, a battery can be used to smooth out
the power produced by a renewable source, or 
consumed by a load such as a computation center.

Like any physical system, batteries degrade over time. The more a battery is cycled,
meaning charged and discharged, the more it ages, and the shorter its
lifetime becomes. The trade-off between short term tasks such as 
arbitrage and renewable smoothing, and a long term task of maximizing battery
lifetime, is the focus of this paper.

\paragraph{Short term objectives.}
We consider two applications, as simple examples of short term objectives.
\begin{enumerate}
\item \textbf{Energy arbitrage.} The battery is charged when electricity is cheap
and discharged when it is expensive. We evaluate the performance of the system
by its revenue.
\item \textbf{Load smoothing.} The battery is used to smooth out the 
electric power demand (\ie, load) consumed by a computation center, so
the net power varies more smoothly over time.
We evaluate the performance of the system using the root
mean square difference between the load and its previous value.
\end{enumerate}
These are just two simple illustrative examples of short term objectives; many others
could be handled by the methods we describe in this paper.

\paragraph{Long term objective.}
The battery has a finite lifetime, which is determined primarily by the number of cycles
it undergoes. We want to maximize the lifetime of the battery, defined as the time
the battery capacity drops below some fraction of its initial value
such as 80\% or 90\%.
Without considering battery aging, the two short term tasks described above 
involve aggressive cycling of the battery,  leading to a shortened lifetime. 
We show how to operate the battery so as to achieve an optimal trade-off
of the short term objective and the long term objective.

\subsection{Related work}
\paragraph{Battery aging models.}
Battery aging is often split into two main effects: calendar aging,
which happens when the battery is resting (no charge or discharge), and
cycle aging, which happens during active charge--discharge 
cycles~\cite{vermeer2021comprehensive}.
Both processes are sensitive to temperature (especially above 30$^\circ$C),
high current rates, and usage patterns like depth of discharge or state of charge 
(SoC)~\cite{xiong2020lithium}. 

When it comes to modeling these effects, three main approaches appear
in the literature. 
Electrochemical models~\cite{keil2020electrochemical, li2020aging, allam2020online} 
try to describe the
internal reactions mathematically (for instance, using the Butler--Volmer
equations~\cite{latz2013thermodynamic}) 
and can be quite detailed but hard to implement in practice.
At the other extreme, empirical models fit observed data to
capture aging trends, but these can fail outside their specific test
conditions and often need large datasets~\cite{pelletier2017battery}. 
Semi-empirical 
models~\cite{serrao2009, suri2016, marano2009lithium, torregrosa2024semi}
blend theoretical ideas with curve-fitting so that the most critical
aging drivers (like temperature, SoC, or C-rate) are handled without
the complexity of a full electrochemical model or the narrow scope of
an entirely data-driven approach.  

Calendar aging mainly depends on how batteries are stored (in terms
of temperature and SoC), with slow chemical reactions gradually eating
away at capacity. Cycle aging, on the other hand, is tied to how
often and how aggressively the battery is charged or discharged. In
actual operation, both processes happen at the same time, so an aging
model typically combines or overlays both~\cite{liu2020evaluation}. 
Semi-empirical models are
especially attractive in this setting because they balance realism
and simplicity, letting us capture key aging behaviors in a way that
can be used in practical control algorithms~\cite{jin2018applicability, miller2022semi}.

\paragraph{Optimization in batteries.}
Current optimization techniques used in battery longevity do not directly focus on 
working with aging models of the batteries. 
Liu \emph{et al.}~\cite{liu_charging_2018} focuses on optimizing 
charging behavior of both CC-CV charging and multi stage CC-CV charging. Nonlinear
optimization techniques are used to determine the optimal current and voltage values
to balance between aging, efficiency loss, and charge time. Chung \emph{et al.}~\cite{chung_optimization_2020} focuses
on reducing the calendar aging of a PEV battery through an optimal charging scheme. 
This paper argues that calendar aging is most important for PEV batteries, 
and proposes a nonlinear interior point method to determine an optimal overnight 
charging scheme for the PEV battery.

Bashir \emph{et al.}~\cite{bashir_lifetime_2017} tackles lifetime maximization of 
lead-acid batteries by formulating key
aging characteristics of the battery as convex formulas. 
Bad recharge, the time since last full
charge, and the lowest state of charge since last recharge are the three primary 
factors in their objective function. A multi-objective convex optimization
problem is then solved to maximize the lifetime of these batteries. 

\paragraph{Model predictive control (MPC).}
Model predictive control goes by several other names, 
such as rolling-horizon planning, receding-horizon control, 
dynamic matrix control, and dynamic linear programming.
Originally developed in the 1960s, MPC offers a
powerful framework for managing constraints on states, inputs, 
and outputs. It has a long history and large literature, and is widely used.
Some early work is~\cite{cutler1979dynamic, garcia1989model}; 
for more recent surveys 
see the papers~\cite{holkar2010overview,mayne2014model,abughalieh2019survey,schwenzer2021review} 
or books~\cite{camacho2007constrained, borrelli2017predictive, grune2017nonlinear, rawlings2017model, rakovic2018handbook}.

Papers describing applications of MPC in specific areas include
HEVs~\cite{huang2017model}
data center cooling~\cite{lazic2018data}, 
building HVAC control~\cite{afram2014theory}, 
wind power systems~\cite{hovgaard2015model}, microgrids~\cite{hu2021model}, 
pandemic management~\cite{carli2020model, peni2020nonlinear}, dynamic hedging~\cite{primbs2009dynamic},
revenue management~\cite{talluri2006theory, bertsimas2003revenue}, 
railway systems~\cite{felez2019model}, aerospace systems~\cite{eren2017model}, 
and agriculture~\cite{ding2018model}. With appropriate forecasting 
(which in many applications is typically simple) 
and choice of cost function, MPC can work well, 
even though it does not explicitly take into account 
uncertainty in the dynamics and cost, or more precisely, 
since it is based on a single forecast of these quantities.

There are many extensions of MPC that attempt to improve performance 
by taking into account uncertainty in the future dynamics and cost. 
Examples include robust MPC~\cite{bemporad2007robust, campo1987robust}, 
min-max MPC~\cite{raimondo2009min}, tube MPC~\cite{mayne2005robust}, 
stochastic MPC~\cite{heirung2018stochastic, mesbah2016stochastic} 
and multi-forecast MPC (MF-MPC)~\cite{shen2021incremental}.

A widely recognized shortcoming of MPC is that it can usually 
only be used in applications with slow dynamics, 
where the sample time is measured in seconds or minutes.
However there exist methods to speed up MPC,
such as computing the entire control law offline 
or using online optimization~\cite{wang2009fast}.

Since  it integrates constraint handling, future forecasting, and feedback adjustment,
MPC is often viewed as a middle ground between exhaustive search methods like
dynamic programming (DP) and simpler real-time strategies such as the
equivalent consumption minimization strategy (ECMS)~\cite{zhang2015comprehensive}. 
This balance between computational tractability and robust performance has made MPC
increasingly popular in both academic and industrial settings for battery
management systems.

\paragraph{Energy arbitrage.}
In economics and finance, arbitrage is the practice of taking advantage of a price difference 
by buying energy from the grid at a low price and selling it back to the grid 
at a higher price~\cite{zafirakis2016value}.
Although it is often assumed to
occur in the day-ahead markets~\cite{staffell2016maximising,wilson2018analysis}, 
arbitrage strategies for
intraday markets have also been considered~\cite{metz2018use}. Multiple studies
assessed how to maximize arbitrage profits~\cite{sioshansi2009estimating}, but the
consistent finding is that the attainable revenues are on their own
insufficient to repay investment in battery storage.
\cite{bradbury2014economic} showed that Li-ion batteries often fail to
surpass 0\% IRR in U.S. markets, though short charge times, lower capital costs, and
ancillary services could enhance returns. \cite{sioshansi2009estimating} found that
arbitrage value in PJM depends on round-trip efficiency, location, and fuel mix,
noting that lower natural gas prices warrant re-evaluation.
\cite{mcconnell2015estimating} highlighted possible profitability in 5-minute dispatch
markets and potential peaker displacement, suggesting that as renewable penetration
depresses prices, arbitrage opportunities may grow in the United States.

The economic value that a battery operator can obtain from arbitrage rests on both
technical and market elements. Among the technical considerations, round-trip
efficiency stands out for its influence on marginal operating costs~\cite{critchlow2017embracing}.
Another crucial factor is the discharge capacity, or energy-to-power ratio, which determines
how much energy a battery can store. Because physical and operational stresses
cause capacity fade over the asset's lifetime~\cite{schmidt2017future}, 
they significantly affect profitability~\cite{he2020economic}. 
Although the replacement cost of a battery is typically incurred
only at the end of its lifespan, Xu \emph{et al.}~\cite{xu2017factoring} notes that aging—driven by
operational choices—should still factor into marginal cost calculations, as it can
alter the operating strategy itself. On the market side, price volatility rather than
average price levels is generally acknowledged as the main determinant of
arbitrage value~\cite{wilson2018analysis}.

\paragraph{Renewable generation and load smoothing.}
A key challenge for large-scale renewable integration is the inherent
fluctuation in power generation, which can cause frequency deviations,
voltage inconsistencies, and high peak loads. 
A common solution is to
pair wind turbines with batteries to 
smooth out these variations, store surplus energy during periods of
high generation, and feed it back when generation drops~\cite{diaz2013energy}.

The issues caused by wind power fluctuations were first discussed in the literature 
in the early 1980s, when commercial wind turbines started being installed more regularly. 
In the first studies, the authors proposed less sophisticated methods 
of power smoothing~\cite{suvire2012improving}. 
In the late 1990s, however, more studies began to consider 
storage systems (mainly fly-wheel and lead-acid batteries) 
to smooth the output power from wind turbines~\cite{jerbi2009fuzzy,el2017second,elkomy2017enhancement}.

In 2009, \cite{khalid2009model,khalid2010model} 
presented a controller design for wind power smoothing purposes 
based on model predictive control. 
They noted that prediction could help improve the economy and security of wind integration into electrical grids.
Thus, a wind power prediction system combined with a battery
was proposed based on measurements from different observation points and communication channels. 
The effectiveness of this approach was assessed through real wind speed data 
from an Australian wind farm comprising 37 wind turbines. 
The results show the capability of the controller to smooth the wind power, 
optimize the maximum ramp rate requirement, and also the state of charge of the battery. 
The study accounted for inefficiencies in batteries in terms of energy conversion but 
did not consider the battery aging.

An analogous problem tackles the challenge of smoothing a load that varies rapidly
over time, such as one that might appear in a data computation center while processing
a job such as training a large language model~\cite{li2024unseen}.
The existing literature explores the use of energy storage that is either integrated
in datacenter uninterruptible power supply (UPS) systems or deployed as
standalone battery banks to enable demand response
services~\cite{urgaonkar2011optimal,govindan2011benefits,govindan2012leveraging,mamun2015physics}. 
Most studies concentrate on minimizing the total
cost of ownership~(TCO)~\cite{wang2014underprovisioning,kontorinis2012managing}, 
defined as the sum of amortized capital
expenditures and operating costs over a prescribed time
horizon~\cite{barroso2019datacenter}. 
Queueing-theoretic Lyapunov optimization has also been used
to derive policies that nearly minimize monthly electricity
bills~\cite{urgaonkar2011optimal,guo2011cutting}. 
While these analyses predominantly assume lead-acid
batteries, a subset of work considers lithium-ion technology. In such demand
response formulations, lithium-ion units are typically modelled as ideal charge
integrators, and aging is captured through depth-of-discharge (DoD) charts or
charge throughput heuristics~\cite{kontorinis2012managing,wang2012energy,aksanli2013distributed,ren2012carbon}.

In 2016, Mamun \emph{et al.}~\cite{mamun2016258}
formulated a multi-objective framework for datacentre demand response
that couples a nonlinear equivalent circuit model with SEI based aging and
tunes feedforward feedback controllers on a lithium-ion battery pack.
Their results reveal an inherent tradeoff between cost savings and battery health;
dead-band PI control mitigates this compromise and remains robust to load
uncertainty and battery pack size.

\subsection{Outline}
We describe the battery model in \S\ref{s-model}
where we use an existing semi-empirical aging model from Suri et al.\cite{serrao2009, suri2016}
and come up with a convex approximation of the aging rate. Next, 
we describe two applications in \S\ref{s-arbitrage} and \S\ref{s-smoothing}
where we use model predictive control for price arbitrage and load smoothing
and give numerical examples.

\section{Cell aging model}\label{s-model}

This section explains battery aging, a mathematical model 
of how a single battery cell ages~\cite{serrao2009,suri2016},
and how we use this model to predict the aging of a larger,
multi-cell battery. To distinguish these two cases we
refer to the former as a ``cell'' and the latter as a ``battery''.

We consider a lithium iron phosphate ($\text{LiFePO}_4$)-graphite battery.
Cells operate by storing chemical potential energy. When a voltage is
applied across the battery terminals, the electrical field causes 
ions to move through the battery's electrolyte, converting the electrical
potential into chemical potential.  This process decays the material of the cell,
especially the cathode and anode. They can crack and oxidize; films can form on
them; ions can become embedded in them. There are complex and detailed
physical models for this process, as predicting lifetime is important
in battery management.

\subsection{Cell charge and current}

We use Suri et al.'s model of lithium-iron phosphate cell aging~\cite{suri2016}.
This model is for a single, 2.5 Ampere-hour (Ah), 3.3 Volt (V) lithium-ion cell.
Physical models such as these model a cell using charging current in
(A) and model capacity and charge in (Ah).
When we model a multi-cell battery, we switch to the more convenient
power and energy units, Watts (W) and Watt-hours (Wh).

We model the cell charging current as constant over time intervals
of length $\delta$ hours, so, \eg, $\delta = 0.25$ means 
$15$ minute intervals.
We denote the time periods as $t=1,2, \ldots$.
We denote the (instantaneous) charge in the cell at the beginning of
interval $t$ as $\tilde q_t$, in units of (Ah).
The charge satisfies $0 \leq \tilde q_t \leq \tilde Q_t$, where $\tilde Q_t>0$
is the capacity in (Ah) of the cell at time in interval $t$.
We assume that $\tilde Q_t$ is known (or measured) at time $t$.
We refer to $\tilde Q_1$ as the initial cell capacity, and $\tilde Q_t$ as the 
cell capacity at time $t$.
The values $\tilde q_t=0$ and $\tilde q_t = \tilde Q_t$ mean that the cell is 
empty and full, respectively.  The empty cell charge $\tilde q_t=0$ refers to the 
lowest charge of the cell over its useful range, and not absolute
zero cell charge.
Similarly, $\tilde q_t = \tilde Q_t$ refers to the largest charge of the cell over its useful
range.

The cell (discharge) current in interval $t$ is denoted $\tilde b_t$, in (A).
Positive values of $\tilde b_t$ correspond to discharging the cell,
and negative values of $\tilde b_t$ correspond to charging the cell.
The cell current $\tilde b_t$ must satisfy $|\tilde b_t| \leq \tilde B$,
where $\tilde B>0$ is the maximum cell charge and discharge current in (A),
given by $\tilde B = \tilde Q_1C$, where $C>0$ is the maximum C-rate of the cell, in inverse
hours (1/h).
The cell dynamics are given by $\tilde q_{t+1} = \tilde q_t - \delta \tilde b_t$. 
Note that $\delta \tilde b_t$ is the total charge, in (Ah), removed from the cell
in interval $t$.

\paragraph{Short and long term quantities.}
We refer to $\tilde b_t$ and $\tilde q_t$ as short term quantities, since they can vary considerably
from interval to interval. 
We refer to $\tilde Q_t$ as a long term (aging) quantity, since it changes very slowly, with 
appreciable change only over a time period measured in months or longer.
In particular, $\tilde Q_t$ can be considered approximately constant over a period on the
order of days.

\subsection{Cell aging}
\paragraph{Loss and loss rate.}
As the cell is used its capacity $\tilde Q_t$ decreases, \ie, 
$\tilde Q_{t+1} \leq \tilde Q_{t}$. 
The lifetime $L$ of the cell in (h) is the time
$L$ when the cell capacity drops below some fixed fraction of its initial value, such
as 90\%,
\ie, $\tilde Q_{L-1} \geq 0.9 Q_1$ and $\tilde Q_L< 0.9 Q_1$. (It is also common
to define lifetime using 80\% of the initial capacity.)
We define the normalized capacity loss as
$l_t = \left( Q_1 - \tilde Q_t \right)/ Q_1$, which we can express as a percentage,
so, \eg, $l_t = 0.07$ means the cell has experienced 7\% capacity loss.
The normalized capacity loss starts at $l_1=0$ (no capacity loss), 
and rises until $l_t > 0.1$, corresponding to cell lifetime.

The cell aging rate is defined as $\rho_t = (l_{t+1} - l_t)/\delta$,
so we have
\[
l_t = \delta \sum_{\tau=1}^{t-1} \rho_\tau.
\]
The aging rate gives the increase in normalized capacity loss per hour, 
and has units (1/h).
A cell with a constant aging rate $\rho$ has a lifetime around 
$L = 0.2/\delta$ hours, which is $2.28 \times 10^{-5}/\delta$
in years, a more commonly used time unit for cell lifetime.
Typical values of $\rho_t$ are on the order of $10^{-6}$ or $10^{-5}$,
corresponding to cell lifetime ranging from around 2 to 20 years.

\paragraph{Loss rate model.}
The aging rate $\rho_t$ depends on how the cell is used, \ie, its
history of charging and discharging up to period $t$.
In this paper we use the semi-empirical aging model 
for a lithium iron phosphate ($\text{LiFePO}_4$)-graphite
cell frequently used in hybrid electric vehicle batteries, given
in~\cite{suri2016}.  (But our methods can be used with 
any other specific cell aging model.)
The model is
\BEQ\label{e-aging-rate}
\rho_t = z \left(\sum_{\tau=1}^t|\tilde b_{\tau}|\delta\right)^{z-1} |\tilde b_t|
\left(\alpha \frac{\tilde q_t}{\tilde Q_t} + \beta \right)
\exp
\left( \frac{-E_a + \eta \frac{|\tilde b_t|}{\tilde Q_t}}{R_g T} \right).
\EEQ
Here $\sum_{\tau=1}^t|\tilde b_{\tau}|\delta$ is 
the accumulated (absolute) charge throughput, up to time $t$,
with units (Ah). 
We note that $\tilde q_t/\tilde Q_t$ is the cell charge normalized to its capacity
(between $0$ and $1$),
and $|\tilde b_t|/\tilde Q_t$ is the instantaneous C-rate of the cell, which is between
$0$ and $C$.

The terms and constants appearing in \eqref{e-aging-rate} are as follows.
\BIT
\item 
The first term models the effect of accumulated charge on aging.
The (unitless) power law exponent is $z = 0.60$.
\item 
The second term models the effect of instantanous cell charge 
of aging rate.
We take $\alpha = 28.966$ and $\beta = 74.112$, with units $(1/\text{Ah}^{z})$.
\item The last term comes from the Arrhenius equation and gives the 
dependence of aging on temperature and instantaneous charging rate.
Here $E_a = 31500$ is the activation energy with units $(\text{J}\text{mol}^{-1})$,
$R_g = 8.314$ is the universal gas constant with units $(\text{J}\text{mol}^{-1}\text{K}^{-1})$,
and $\eta = 152.500$, with
units $(\text{J}\text{mol}^{-1}\text{h})$.
The absolute temperature $T$ is given in degrees Kelvin (K).
\EIT

The accumulated charge and cell capacity $\tilde Q_t$ are long term quantities,
which do not change much over periods of a few days (after the initial 
few months, in the case of the accumulated charge).
Combining the long term quantities we can
rewrite the capacity loss rate~\eqref{e-aging-rate} as a function
of the short term quantities as
\BEQ\label{e-aging-rate-simplified}
\rho_t = \mu_t |\tilde b_t| (1+\nu_t \tilde q_t) \exp \lambda_t |\tilde b_t|,
\EEQ
where 
\[
\mu_t = 
\beta \exp \left( \frac{-E_a}{R_gT} \right)
z \left(\sum_{\tau=1}^t|\tilde b_{\tau}|\delta\right)^{z-1},
\qquad \nu_t = \frac{\alpha}{\beta \tilde Q_t}, \qquad
\lambda_t = \frac{\eta}{R_g T \tilde Q_t}
\]
are long term, slowly varying quantities, known at time $t$.
Over short time periods (\eg, a few days), the coefficients
$\mu_t$, $\nu_t$, and $\lambda_t$ can be considered constant, so the 
aging loss rate depends only on $|\tilde b_t|$ and $\tilde q_t$.
It increases with $|\tilde b_t|$ proportional to $|\tilde b_t|\exp \lambda_t |\tilde b_t|$,
and proportional to $1+\nu_t \tilde q_t$ with $\tilde q_t$.

\paragraph{Short term approximation.}
We develop here a convex approximation of $\rho_t$ in \eqref{e-aging-rate-simplified} 
that can be used in the short term.
Taking the first order Taylor expansion of \eqref{e-aging-rate-simplified} with
respect to $|\tilde b_t|$ and $\tilde q_t$, around 
the point $|\tilde b_t|=0$ and $\tilde q_t = \tilde Q_t/2$, we 
obtain
\BEQ\label{e-aging-rate-approx}
\hat \rho_t = \mu_t \left(1 + \nu_t \frac{\tilde Q_t}{2} \right) |\tilde b_t|.
\EEQ
This approximation seems crude, but predicts aging rate reasonably accurately.

\subsection{Battery aging model}
The above model is for a single battery cell. Larger batteries are made up of
many cells.  We consider a multi-cell battery consisting of $N$ cells,
either passively wired in a series-parallel arrangement, or with
active battery power control, or any combination.
We model the multi-cell battery with capacity given in energy (Wh) and
charge/discharge given in power (W).
We assume that the battery cells work in a balanced way, \ie, all
have the same current, charge, and capacity at each time.

We denote the cell charging current as $\tilde b_t$ in (A), 
cell charge as $\tilde q_t$ in (Ah), and cell capacity as $\tilde Q_t$
in (Ah), as above. We denote the battery charging power in (W) as $b_t$,
the battery charge as $q_t$ in (Wh), and the battery
capacity as $Q_t$ in (Wh).
These are related as
\[
b_t = 3.3 N \tilde b_t , \quad
q_t = 3.3 N \tilde q_t , \quad
Q_t = 3.3 N \tilde Q_t,
\]
where $3.3$ (V) is the cell voltage.
(Here we make the reasonable approximation that the cell voltage is constant over its useful
charge range.)
The maximum C-rate of the battery is the same as that of a single cell.

\subsection{Numerical example}
\begin{figure}
\centerline{\includegraphics[width=0.8\columnwidth,keepaspectratio,clip=true]{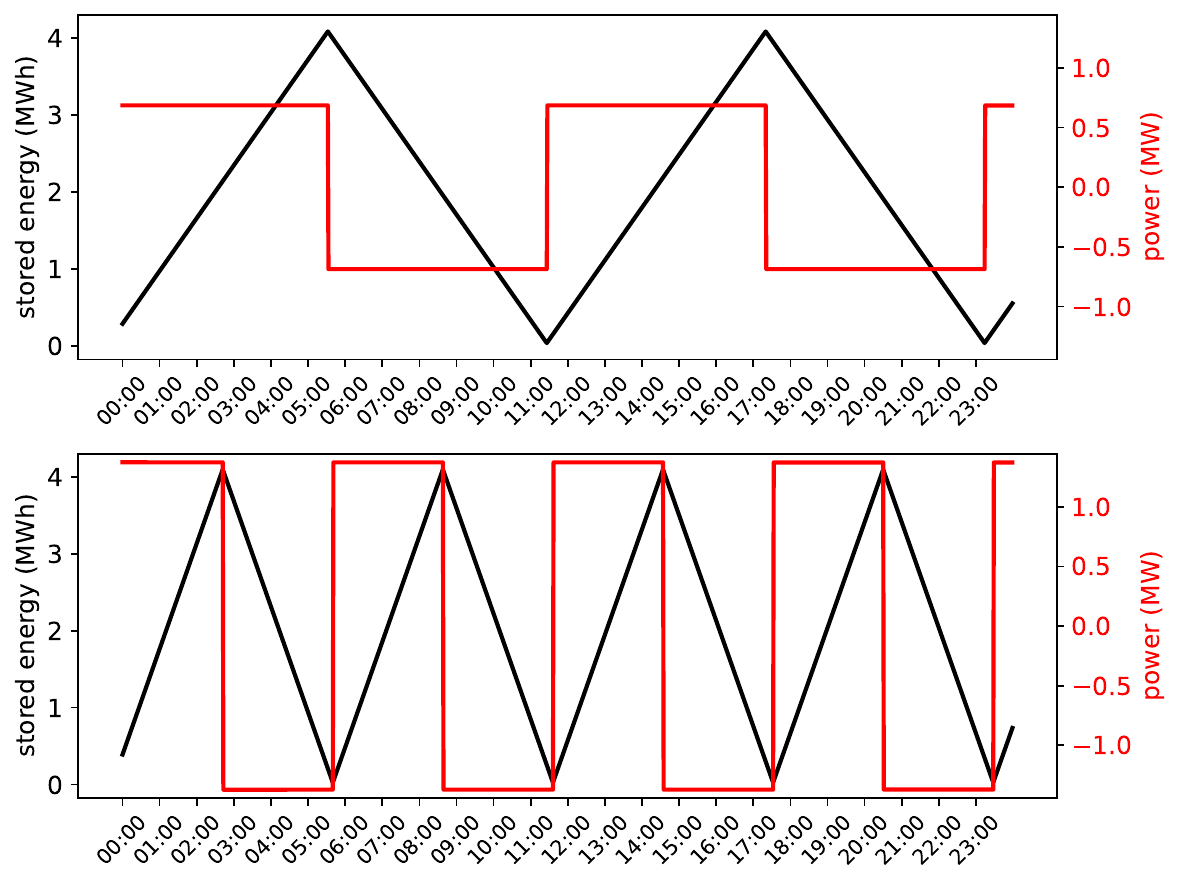}}
\caption{Charging profiles. \emph{Top.} 2 cycles per day. \emph{Bottom.} 4 cycles per day.}
\label{f-charging_profiles}
\end{figure}

\begin{figure}
\centerline{\includegraphics[width=0.8\columnwidth,keepaspectratio,clip=true]{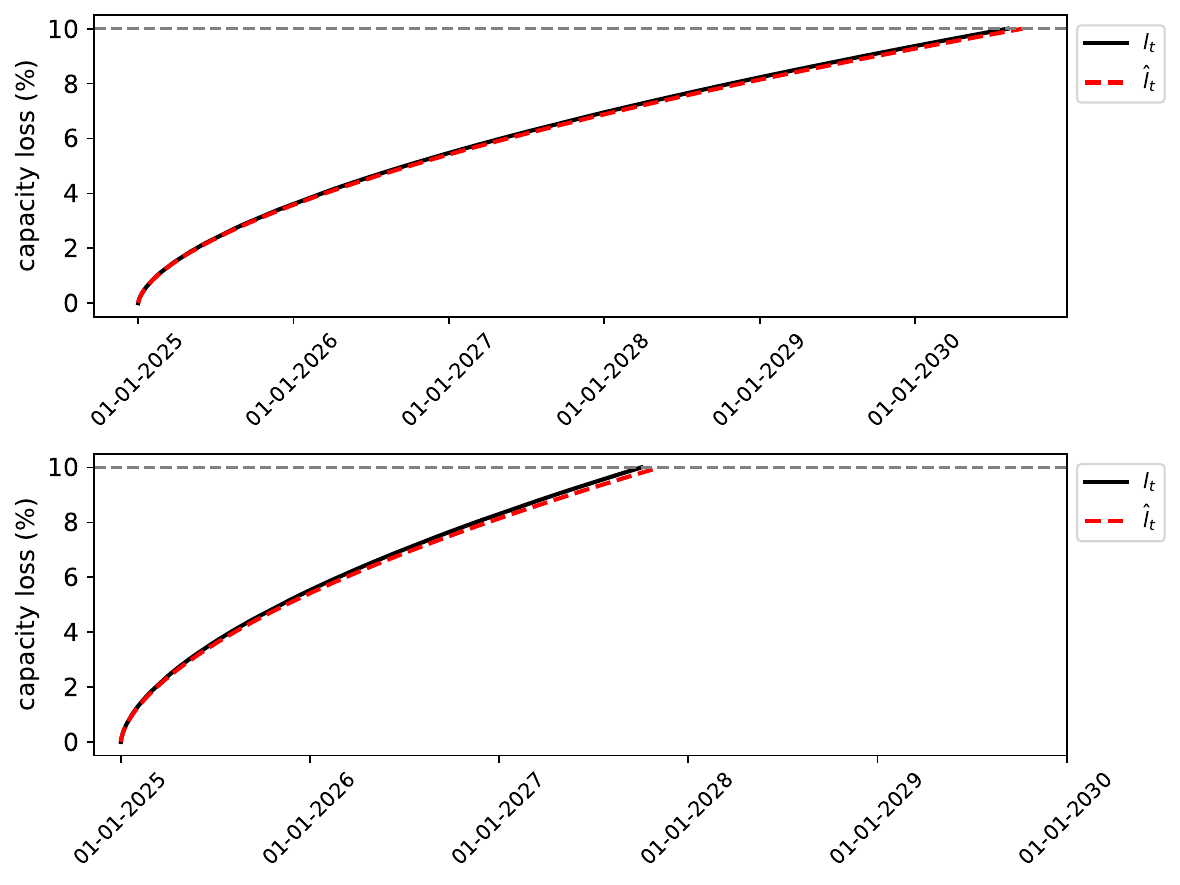}}
\caption{Actual and approximate aging. \emph{Top.} 2 cycles per day. \emph{Bottom.} 4 cycles per day.}
\label{f-loss_profiles}
\end{figure}

In this section we demonstrate the aging model \eqref{e-aging-rate-simplified} 
and the short-term approximation \eqref{e-aging-rate-approx} with examples.
The battery has $N=5000000$ cells, corresponding to a capacity $Q_1 = 4.125$ (MWh).
The maximum C-rate is set to $0.33$ (1/h) which means the battery can completely
charge or discharge in $3$ hours.
We consider two simple charging profiles, shown in figure~\ref{f-charging_profiles}.
The first one fully charges and discharges the battery, with constant 
charge/discharge current, with $2$ cycles per day.
The second one does the same, with $4$ cycles per day.
The associated capacity losses are shown in figure~\ref{f-loss_profiles}.
In these plots we show the aging capacity loss calculated using the 
exact aging model \eqref{e-aging-rate-simplified} and the approximate aging model
\eqref{e-aging-rate-approx}.
The approximate aging model matches the exact aging model well: 
Table~\ref{table:lifetime_vs_charging}
shows it predicts lifetimes within a few percent of the more accurate model.
\begin{table}
\centering
\begin{tabular}{r|r|r|}
\makecell[l]{Charging profile\\ (cycles/day)} &
\makecell[r]{Actual lifetime\\ (years)} &
\makecell[r]{Approximate lifetime\\ (years)} \\ 
\hline
2 & 5.60 & 5.70 \\
4 & 2.75  & 2.85  \\
\end{tabular}
\caption{Lifetime for different charging profiles, calculated using
exact and the approximate aging models.}
\label{table:lifetime_vs_charging}
\end{table}

\clearpage
\section{Arbitrage}\label{s-arbitrage}
\subsection{Problem setup}\label{ss-arbitrage-setup}
In our first example the short term task is arbitrage with time-varying energy prices.
The (nonnegative) price of electricity in period $t$ is $p_t$ in (USD/MWh).
The battery is connected to the grid, and
a payment of $p_tb_t\delta$ is received in period $t$; negative payments are amounts 
we pay to the grid operator.
The objective is to choose $b_t$ to maximize the average (net) payment over a day 
given by $\frac{1}{T^{\text{day}}} \sum_{\tau=1}^{T^{\text{day}}} p_\tau b_\tau \delta$,
where $T^{\text{day}}$ is the number of time periods in a day 
and $\delta$ is the period length in hours (h).
The prices are known far enough ahead of time (\eg, on day), that we can 
consider them known.

Without considering battery aging, the greedy strategy is to discharge the battery 
as much as possible when the price is high, 
and charge it as much as possible when the price is low.
In that case, we are only subject to physical limits such as maximum discharge rate, 
battery capacity, current storage, etc.
However this strategy will lead to aggressive cycling of the battery and will shorten its lifetime.
That is why we also consider the total revenue over the lifetime of the battery,
$\sum_{\tau=1}^{T^{\text{EOL}}} p_\tau b_\tau \delta$, where $T^{\text{EOL}}$ is the 
end of life of the battery.

\subsection{Short term MPC method}\label{ss-arbitrage-mpc}
At time $t$ we are given
battery storage $q_t$, battery capacity $Q_t$, 
and compute the approximate aging rate coefficient $\mu_t(1 + \nu_t \frac{Q_t}{2})$
using the battery model in \S\ref{s-model}.
We consider a horizon of $H$ periods, and assume that together with $Q_t$, 
the aging rate coefficient is constant over this horizon.
We assume that the price of electricity $p_t$ is known for the next $H$ periods
and that they are nonnegative.
To find the optimal battery discharge $b_{\tau}$, 
and storage $q_{\tau}$
for $\tau=t+1,\ldots,t+H$, we solve the problem
\[
\begin{array}{ll}
\text{minimize} 	& \frac{1}{H} \sum_{\tau=t+1}^{t+H} \left(-p_{\tau} b_{\tau} \delta
					+ \gamma \mu_{t} (1 + \nu_t \frac{Q_t}{2}) |b_{\tau}| \right) + \eta \left(q_{t+H} - \frac{Q_t}{2} \right)^2\\
\text{subject to} 	& q_{\tau} = q_{\tau-1} - \delta b_{\tau-1}, \quad \tau=t+1, \ldots, t+H \\
					& |b_{\tau}| \leq C Q_t, \quad \tau = t+1, \ldots, t+H \\
					& 0 \leq q_{\tau} \leq Q_t, \quad \tau = t+1, \ldots, t+H	
\end{array}
\]
where $\gamma >0$ is a trade-off parameter between revenue and battery aging,
and $\eta>0$ is a parameter that penalizes deviation of the terminal battery charge from 
$Q_t/2$, half capacity.
Solving this optimization problem gives us a plan for operating the battery over
the next $H$ periods. Our policy uses the first battery power $b_t$ in our plan.

\subsection{Data}\label{ss-arbitrage-data}
We use hourly local marginal prices (LMP) data
for the day-ahead market in the Electric Reliability Council of Texas (ERCOT)
North Hub for the year 2012. The data is available at~\cite{ercotprice2012}.
We use data from January 1, 2012 to December 31, 2012 and have $8784$ data points.
We run simulations until the battery end of life, which ranges from $6$ to $22$ years.

To simulate multiple years, we repeat the data from 2012, removing the data
from February 29th for non-leap years.
Repeating the same data for multiple years removes yearly variations 
in the data, which varies considerably from year to year.
For instance, in February 2021, the price of electricity in Texas spiked to
$\sim9000$ (USD/MWh), two orders of magnitude higher than the average price in 2012, 
due to a winter storm. To reduce the effect of such outliers and 
yearly variations, we use the same data for multiple years.

\begin{figure}
\centerline{\includegraphics[width=0.8\columnwidth,keepaspectratio,clip=true]{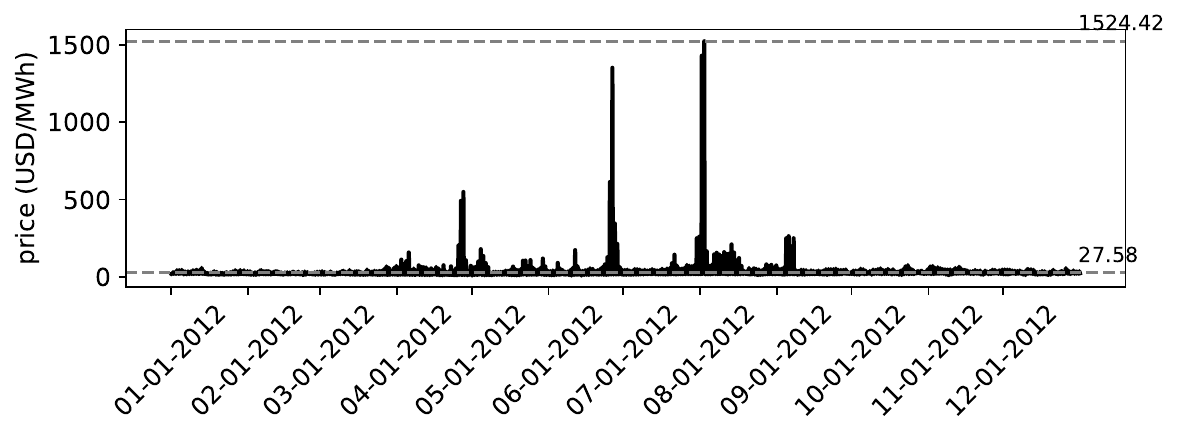}}
\caption{ERCOT hourly local marginal prices for day-ahead market in 2012. 
The dashed lines show the mean and the maximum.}
\label{f-yearly_price_data}
\end{figure}

To visualize seasonal variation, we first look at prices at the yearly scale,
shown in figure~\ref{f-yearly_price_data}.
The mean is $27.58$ (USD/MWh) and its standard deviation is $39.36$ (USD/MWh).
The maximum price is $1524.42$ (USD/MWh).
We observe that there are $2$ major peaks, both in summer months.

\begin{figure}
\centerline{\includegraphics[width=0.8\columnwidth,keepaspectratio,clip=true]{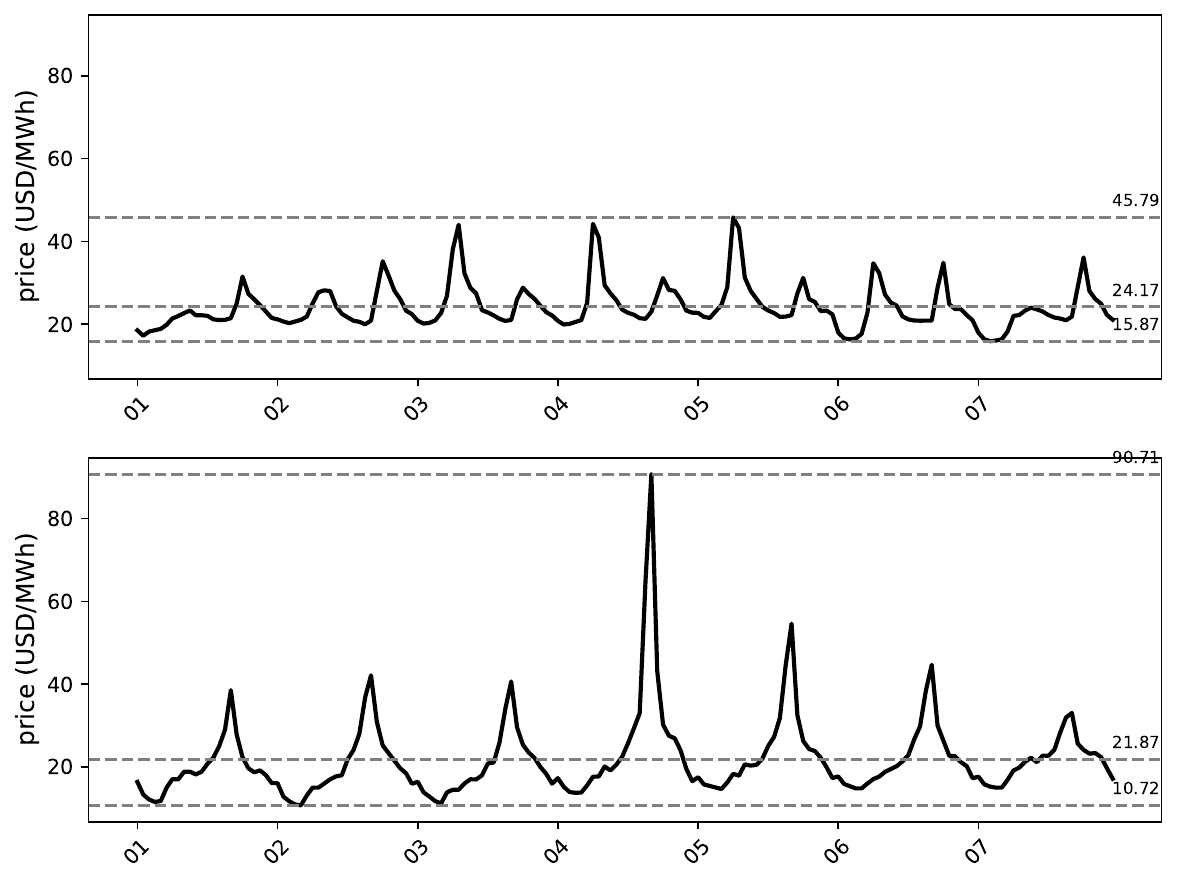}}
\caption{Two different weeks of price data. \emph{Top.} January 2012.
\emph{Bottom.} June 2012.}
\label{f-weekly_price_data}
\end{figure}

Next, we look at the data at a finer scale.
Figure~\ref{f-weekly_price_data} shows the prices for the first
week of January $2012$ and the first week of June $2012$.
We observe several expected phenomena, \eg,
prices are higher during the day than at night, and a bit higher in winter
than in summer. We can also see some small variation over a week.
One interesting observation is that the shape of the daily demand
in winter differs considerably from the shape of the daily demand
in summer. In winter we see a double bump, with peaks in the morning
and afternoon, while in summer we see a smoother daily variation
with one peak in the early afternoon. 
We interpret this as the effect of heating in winter, which causes
a peak in the morning, and air conditioning in summer, which causes
a peak in the afternoon.
Also, even if the daily demand
curves are similar for the same season, they show some variation.
For instance, the daily demand on 06-04-2012 is higher with a sharper peak
around noon compared to the other days of the week.

We assume that we know the prices for the next $H$ periods
beforehand regardless of the time of day, so we do not need to rely on 
a forecasting model to predict the prices.
However, in practice, day-ahead market 
prices in ERCOT are only available at around 1:00 PM for 00:00 AM to 11:00 PM on the next day.
As a result, we would need to use a forecaster for at least some of the hours.

\subsection{Simulation results}\label{ss-arbitrage-results}
Using the MPC method in \S\ref{ss-arbitrage-mpc} and data in \S\ref{ss-arbitrage-data},
we change the cost of battery aging $\gamma$ and simulate hourly (\ie, \ $\delta = 1$ (h))
until the end of life of the battery.  We take horizon $H=24$, \ie, one day.
We found that the results are not sensitive to
the value of $\eta$, and use $\eta=1$.

\begin{figure}
\centerline{\includegraphics[width=0.8\columnwidth,keepaspectratio,clip=true]{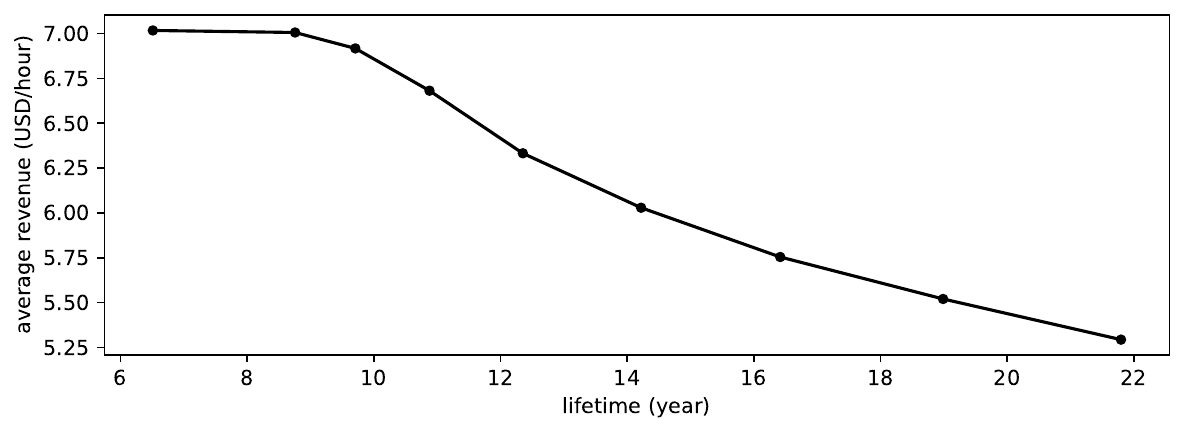}}
\caption{Average hourly revenue versus battery lifetime.}
\label{f-total_and_average_revenue}
\end{figure}

The total revenue and average revenue versus battery lifetime
is shown in figure~\ref{f-total_and_average_revenue}.
We observe that the total revenue increases with increasing battery lifetime
while the average revenue decreases.
This is expected, since with a shorter battery lifetime, the battery is 
allowed to cycle more aggressively, which increases the hourly revenue.
However, this comes at the cost of a shorter battery lifetime and lower
total revenue.
Hence there is a trade-off between the total revenue and the average revenue.

Since total revenue is a monotonically increasing function of battery lifetime, and 
we would get higher revenue with a longer battery lifetime, we focus on the 
net present value (NPV) of the revenue over the lifetime of the battery.
The NPV is a discounted sum of the revenue over the lifetime of the battery, 
where the interest rate corresponds to the cost of capital.
For a given interest rate $i$, NPV is given by
$\text{NPV}(i) = \sum_{\tau=1}^{T^{\text{EOL}}} \frac{p_\tau b_\tau \delta}{(1+i)^\tau}$.

\begin{figure}
\centerline{\includegraphics[width=0.8\columnwidth,keepaspectratio,clip=true]{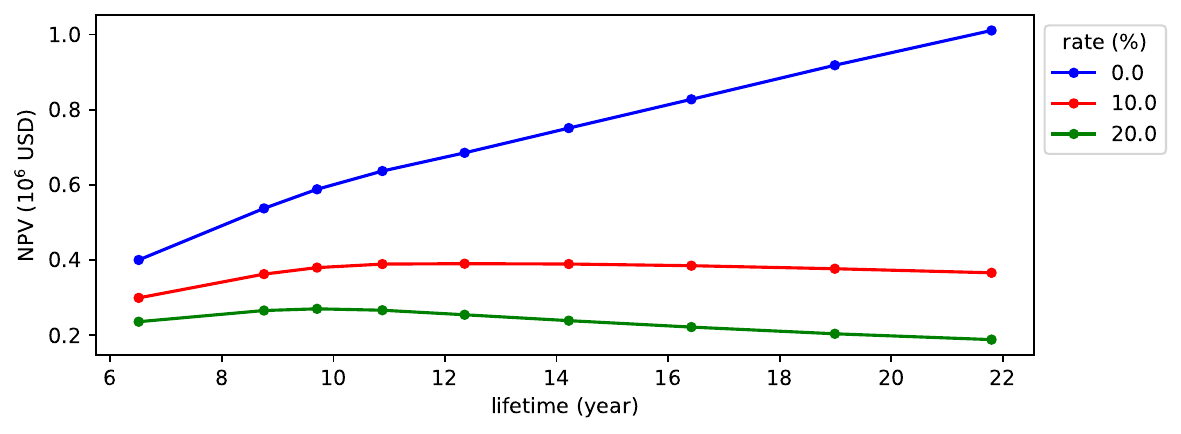}}
\caption{Net present value of revenue versus battery lifetime.}
\label{f-npv_vs_lifetime}
\end{figure}

We plot the NPV versus battery lifetime in figure~\ref{f-npv_vs_lifetime}
and show for interest rates $i=0\%,10\%,20\%$.
We observe that with $20\%$ interest rate, the NPV is maximized with a battery
lifetime of around $10$ years.

\begin{figure}
\centerline{\includegraphics[width=0.8\columnwidth,keepaspectratio,clip=true]{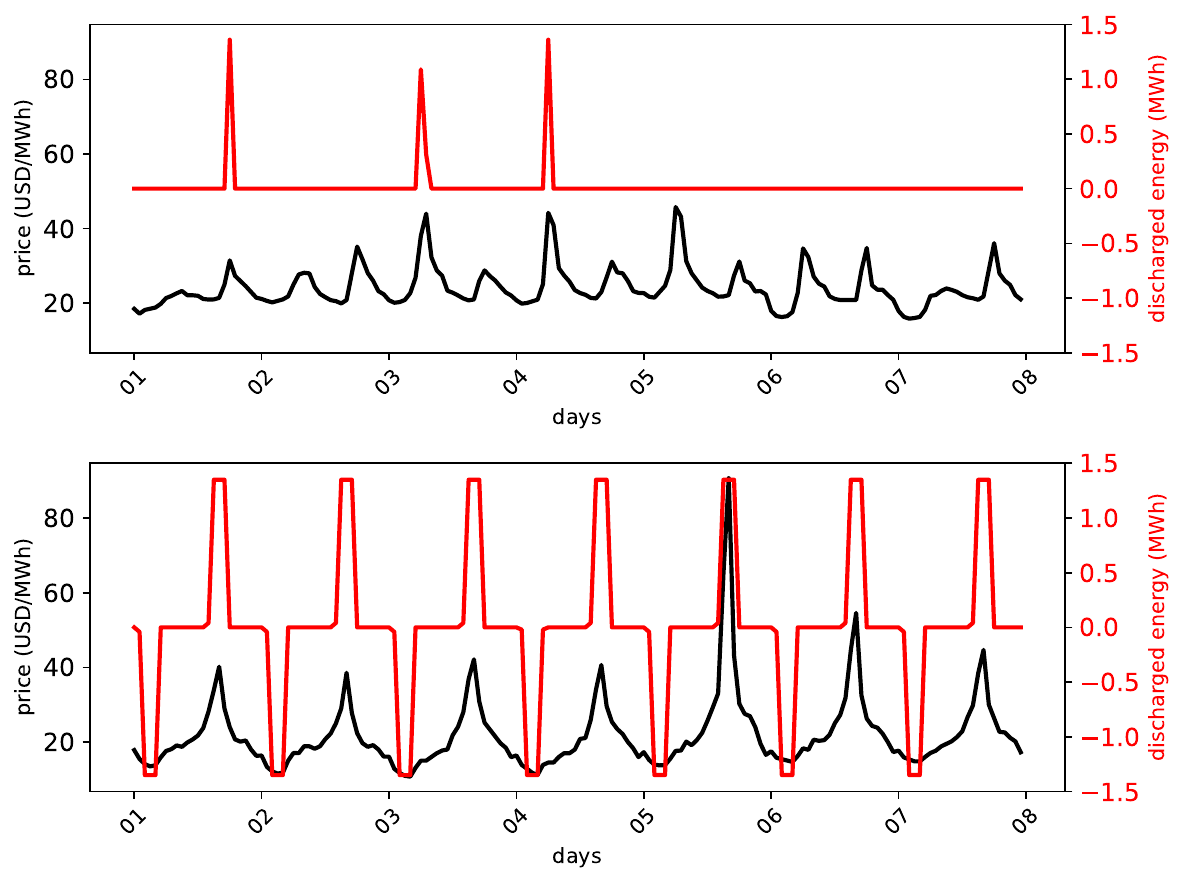}}
\caption{Weekly discharge profile for a battery with $8.75$ years lifetime. \emph{Top.} January 2012.
\emph{Bottom.} June 2012.}
\label{f-weekly_discharge_profile_arbitrage}
\end{figure}

Next, we focus on discharge profiles at a smaller timeframe and 
look at the weekly discharge profile for a battery with an $8.75$ year lifetime
in figure~\ref{f-weekly_discharge_profile_arbitrage}.
As expected, we observe that discharge periods coincide with high price periods,
and charge periods coincide with low price periods.
In fact, during the summer the battery is discharged almost every day around noon,
where the price is highest, and charged in the early morning, where the price is lowest.
In winter, the battery is idle most of the time, and is discharged/charged
only a few times a week.

\clearpage
\section{Load smoothing}\label{s-smoothing}
\subsection{Problem setup}\label{ss-smoothing-setup}
In our second example the short term task is to smooth out a time-varying load.
We have a load $w_t$ in Watts (W),
and we operate the battery in parallel, so the total power is $z_t = w_t + b_t$.
We want $z_t$ to be smooth, as judged by the RMS difference
\[
\mathcal D = \left( \frac{1}{T-1} \sum_{t=1}^{T-1} (z_{t+1}-z_t)^2 \right)^{1/2},
\]
with smaller values better.
At time period $t$, the load $w_t$ is known; future values $z_{t+1}, z_{t+2}, \ldots$
are not known, but can be forecasted.

\subsection{Short term MPC method}\label{ss-smoothing-mpc}
We use an MPC policy.
At time $t$, we are given battery storage $q_t$,
battery capacity $Q_t$, load $w_t$, and aging parameter 
$\mu_t\left(1+\nu_t \frac{Q_t}{2}\right)$. We assume the battery capacity and 
aging parameter are constant over our short term horizon $H$.
We also have load forecasts $\hat{w}_{\tau \mid t}$ for $\tau=t+1,\ldots,t+H$.
We will discuss how we obtain these forecasts
later in \S\ref{ss-load-forecasts}.
We solve the short term planning problem
\[
\begin{array}{ll}
\text{minimize} 	& \frac{1}{H} \sum_{\tau=t+1}^{t+H} \left((z_{\tau}-z_{\tau-1})^2 
					+ \gamma \mu_{t} (1 + \nu_t \frac{Q_t}{2}) |b_{\tau}| \right) + \eta \left(q_{t+H} - \frac{Q_t}{2} \right)^2 \\
\text{subject to} 	& z_{\tau} = \hat{w}_{\tau \mid t} + b_{\tau}, \quad \tau=t+1, \ldots, t+H \\
					& z_{t} = w_t + b_t \\
					& 0 \leq z_{\tau}, \quad \tau=t+1, \ldots, t+H \\ 
					& q_{\tau} = q_{\tau-1} - \delta b_{\tau-1}, \quad \tau=t+1, \ldots, t+H \\
					& |b_{\tau}| \leq C Q_t, \quad \tau=t+1, \ldots, t+H \\
					& 0 \leq q_{\tau} \leq Q_t, \quad \tau=t+1, \ldots, t+H
\end{array}
\]
where $\gamma > 0$ is a trade-off parameter between smoothing and battery aging,
and $\eta>0$ penalizes deviation of the final battery storage from half capacity.
The solution of this problem gives us a plan for $b_t, \ldots, b_{t+H}$; we use as $b_t$
the first battery charge in this plan.

\subsection{Data}\label{ss-smoothing-data}
\begin{figure}
\centerline{\includegraphics[width=0.8\columnwidth,keepaspectratio,clip=true]{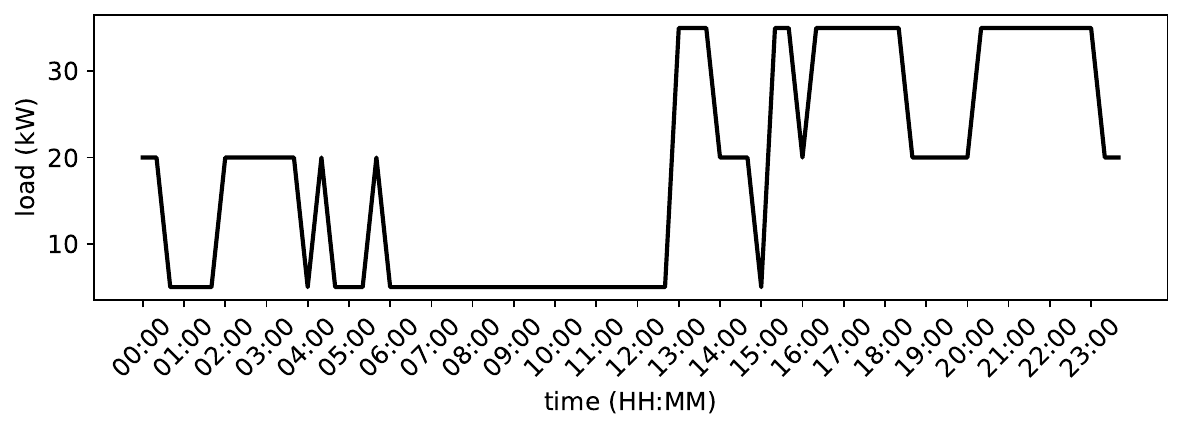}}
\caption{\parbox{0.8\columnwidth}{Simulated load data for January 1st 2020.}}
\label{f-daily_load_data}
\end{figure}

We use a battery with $N=15000$ cells, corresponding to $123.75$ (kWh) initial battery capacity.
Given that most commercial residential batteries have a capacity of $13.5$ (kWh),
this is equivalent to having $9$ batteries. Maximum charge/discharge rate is $C=0.3$ (1/h) 
which implies that the battery can completely charge/discharge in around $3$ hours.
We use simulated data modeled after the 
consumption of a large language model (LLM) training job
in the MIT Supercloud Dataset as described in Li \emph{et al.}~\cite[Figure 6]{li2024unseen}.
We use $20$ minute periods, and three discrete values of 
load: $5$ (kW), $20$ (kW), and $35$ (kW), corresponding to different computation states.
We refer to these as the low, medium, and high power states, respectively.
We use a Markov model for the computation states, taken (roughly) from
\cite[Figure 6]{li2024unseen}, with transition matrix
\[
P = 
	\begin{bmatrix}
	0.79 & 0.22 & 0.00 \\
	0.05 & 0.72 & 0.40 \\
	0.16 & 0.06 & 0.60
	\end{bmatrix}
\]
where $P_{ij}$ is the probability of transitioning from state~$j$ to state~$i$.
We generate $25$ years of data starting from 01-01-2018 00:00 PST,
which gives $657000$ data points.
The asymptotic state probabilities are $0.40$, $0.38$, and $0.22$.  The asymptotic
average load power is $17.25$ (kW).
Figure~\ref{f-daily_load_data} shows the simulated load data for 
January 1st 2020.

\subsection{Load forecasts}\label{ss-load-forecasts}
We use a simple conditional mean forecast $\hat{w}_{\tau \mid t}$ 
for $\tau=t+1,\ldots,t+H$,
obtained as 
\[
\hat w_{\tau \mid t} = (5,20,35) P^{t-\tau} s_t, \quad \tau=t+1, \ldots, t+H,
\]
where $s_t$ is the current state, \ie, $s_t = (1,0,0)$ in the low power
state, $s_t = (0,1,0)$ is the medium power state, and $s_t= (0,0,1)$ in
the high power state.

\begin{figure}
\centerline{\includegraphics[width=0.8\columnwidth,keepaspectratio,clip=true]{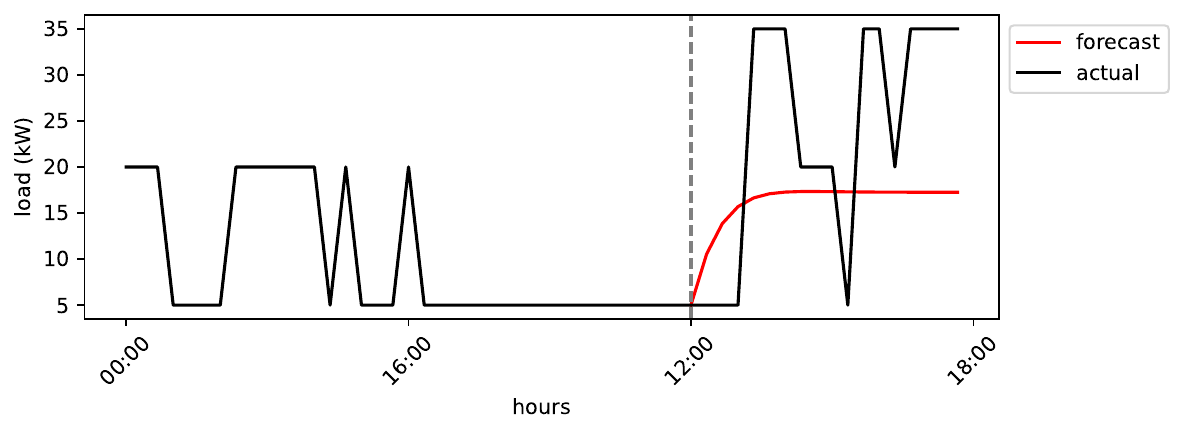}}
\caption{\parbox{0.8\columnwidth}{Actual load and forecasts on January 1st 2020.}}
\label{f-daily_load_forecasts}
\end{figure}

We show the load data and $6$ hour ahead forecasts for January 1st 2020
in figure~\ref{f-daily_load_forecasts}.
We observe that the forecasts converge to 
steady state after around $2$ hours.
\subsection{Simulation results}\label{ss-smoothing-results}
Using the MPC method in \S\ref{ss-smoothing-mpc}, data in \S\ref{ss-smoothing-data}
and forecast in \S\ref{ss-load-forecasts},
we change the cost of battery aging $\gamma$ and 
simulate in $20$ minute intervals (\ie, \ $\delta = 1/3$ (h))
until the end of life of the battery.
Horizon $H$ is set to $6$ hours. 
We observed that the choice of $\eta$ did not affect the results significantly
and we set it to $\eta = 0.5$ for all simulations.

\begin{figure}
\centerline{\includegraphics[width=0.8\columnwidth,keepaspectratio,clip=true]{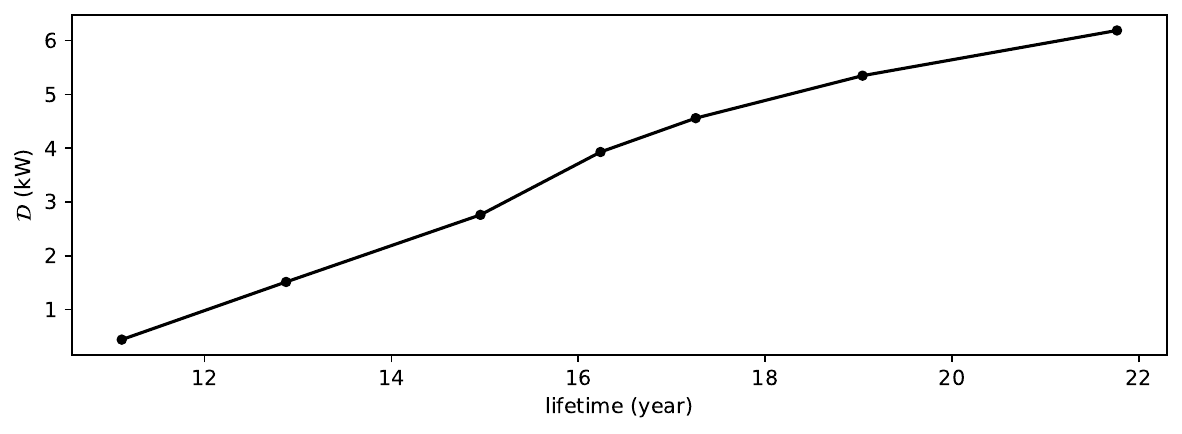}}
\caption{$\mathcal D$ versus battery lifetime.}
\label{f-smoothing_tradeoff_curve}
\end{figure}

RMSD (\ie, $\mathcal D$) of the smoothed load versus battery lifetime
is shown in figure~\ref{f-smoothing_tradeoff_curve}.
$\mathcal D = 10.29$ (kW) for $25$ years of simulated load.
We observe that with a battery lifetime of $11$ years, 
$\mathcal D = 0.44$ (kW) and load is smoothed out approximately $23$ times.
On the other hand, with a battery lifetime of $22$ years,
$\mathcal D = 6.20$ (kW) and load is smoothed out approximately $1.7$ times.
Hence $\mathcal D$ increases with increasing battery lifetime.
This is expected, since with a shorter battery lifetime, the battery is 
allowed to cycle more aggressively, which results in a smoother signal
and lower $\mathcal D$. However, this comes at the cost of a shorter battery lifetime.

\begin{figure}
\centerline{\includegraphics[width=0.8\columnwidth,keepaspectratio,clip=true]{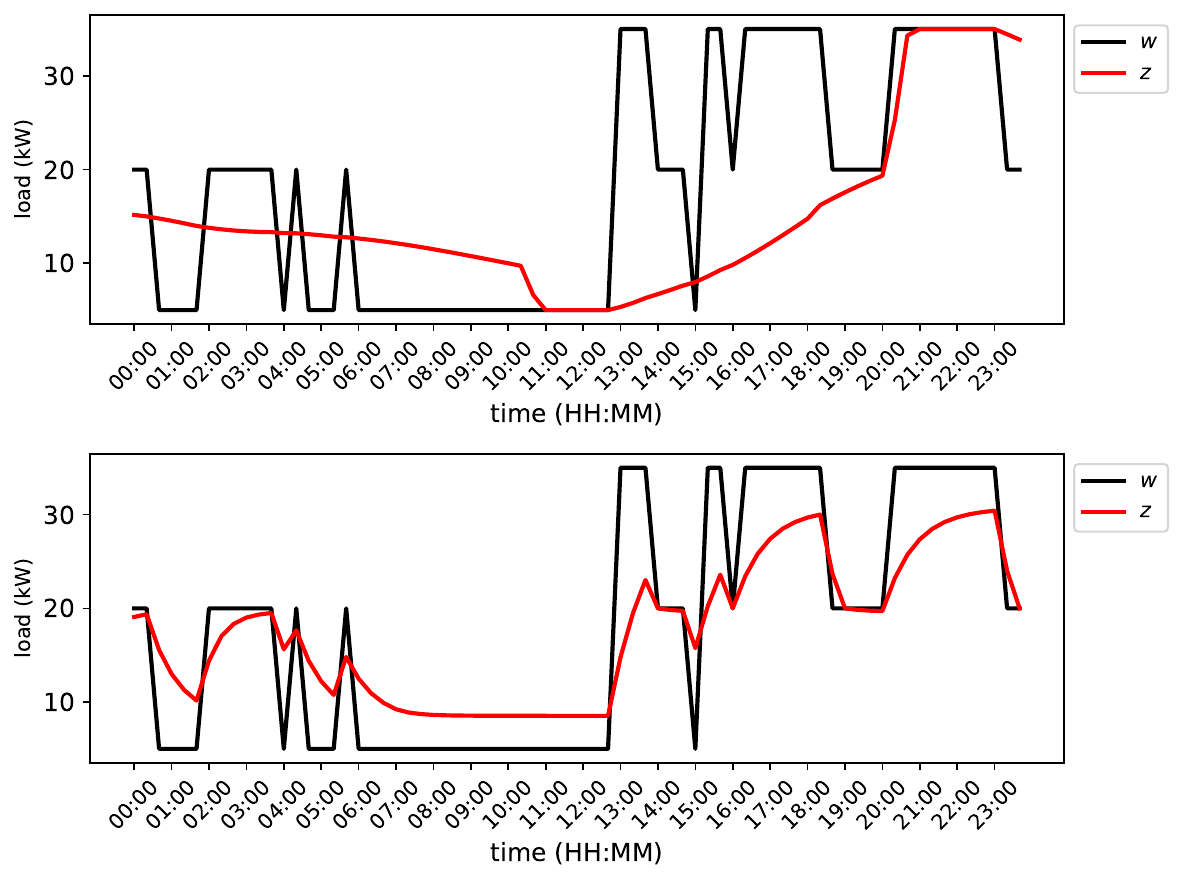}}
\caption{Daily smoothing profile using forecasts. \emph{Top.} $11$ years battery lifetime.
\emph{Bottom.} $15$ years battery lifetime.}
\label{f-daily_smoothing_profile}
\end{figure}

We focus on smoothing profiles at a smaller timeframe and 
look at the daily load smoothing profile for a battery with an $11$ year lifetime
and $15$ year lifetime in figure~\ref{f-daily_smoothing_profile}.
The smoothed load changes more slowly than the original load, 
and as the battery lifetime increases, the smoothing decreases.

\clearpage
\section{Conclusion}\label{s-conclusion}
In this paper we formulated the battery management problem considering two
competing objectives: 
short term goals (with specific examples energy arbitrage and load smoothing)
and long term battery lifetime maximization. 
We adopted an existing
semi-empirical battery aging model for lithium iron phosphate cells and
developed a convex optimization-based control strategy 
using a simplified yet accurate convex approximation of the aging rate.

Our approach employs MPC, 
leveraging known future prices or load forecasts to optimally balance 
short term task performance and battery longevity. 
Through extensive numerical simulations, we clearly
demonstrated the trade-off between aggressive short term battery cycling and
battery lifetime. 
We observed that aggressive cycling significantly increased
hourly revenues and reduced load fluctuations, 
but at the cost of a shorter battery lifetime. 
Conversely, conservative strategies prolonged battery life
but yielded lower short term performance gains.

Key novelties of our work include the convex approximation of the battery aging
model, enabling computationally efficient optimization within the MPC
framework, and systematic quantification of the lifetime-performance trade-off
under realistic conditions. 
Our results underscore the importance of
aging-aware optimization strategies, providing a practical and scalable
solution that bridges the gap between detailed battery aging models and
real-time battery management.
Open source implementation and data used in this work are available at~\cite{cvxbattery2025}.
 
\section*{Acknowledgments}
We would like to thank Prof.~Simona Onori 
for discussions of the battery aging model,
Dr.~Eric Sager Luxenberg for discussion of 
an initial formulation and Maximilian Schaller 
for careful reading and corrections of the manuscript
and code. 

\clearpage
\bibliography{references}
\clearpage

\end{document}